\documentclass[11pt]{article}

\author{}
\date{}

\usepackage{bbm,amsmath,amsthm,amsfonts,amssymb}

\usepackage{graphicx}

\newcommand{\captionfonts}{\footnotesize}
\makeatletter  
\long\def\@makecaption#1#2{%
  \vskip\abovecaptionskip
  \sbox\@tempboxa{{\captionfonts #1: #2}}%
  \ifdim \wd\@tempboxa >\hsize
    {\captionfonts #1: #2\par}
  \else
    \hbox to\hsize{\hfil\box\@tempboxa\hfil}%
  \fi
  \vskip\belowcaptionskip}
\makeatother   

\newcommand{\nwc}{\newcommand}

\newtheorem{prop}{Proposition}[section]

\newtheorem{theorem}[prop]{Theorem}

\nwc{\R}{\mathbb R}
\nwc{\Z}{\mathbb Z}
\nwc{\N}{\mathbb N}

\newcommand{\ignore}[1]{}

\nwc{\eps}{\varepsilon}
\nwc{\re}{Re\,}

\nwc{\wto}{\rightharpoonup}

\nwc{\ds}{\displaystyle}
\newcommand {\bedis} {\begin{displaymath}}
\newcommand {\edis} {\end{displaymath}}
\newcommand{\newbeqna} {\renewcommand {\arraystretch} {2}
                        \begin {displaymath} \begin {array}{crcl}}
\newcommand{\neweqna}{\end{array} \end {displaymath}}

\newcommand{\fbeqna}{\renewcommand {\arraystretch} {1.3}
\begin {displaymath}\begin{array}{rcll}}
\newcommand{\feqna}{\end{array}\end{displaymath}}

\newcommand {\beqna} {\begin{eqnarray*}}
\newcommand {\eqna} {\end{eqnarray*}}
\newcommand {\beqn} {\begin{eqnarray}}
\newcommand {\eqn} {\end{eqnarray}}

\begin{document}
\title{
Self-similar solutions with fat tails for a coagulation equation with
diagonal kernel
}

\author{%
B. Niethammer%
\footnote{Mathematical Institute, University of Oxford, 24-29 St.
Giles, Oxford, OX1 3LB, England }
and
J. J. L. Vel\'{a}zquez\footnote{ICMAT (SCIC-UAM-UC3M-UCM), 
 Facultad de Ciencias Matem\'{a}ticas, Universidad Complutense,
Madrid 28040, Spain}.
}%
\maketitle

\begin{abstract}

We consider  self-similar solutions of Smoluchowski's coagulation
equation with a diagonal kernel of homogeneity $\gamma < 1$.
We show that there exists a family of 
second-kind self-similar solutions with power-law behavior
$x^{-(1+\rho)}$ as $x \to \infty$ with $\rho \in (\gamma,1)$. To our knowledge this is the first example of a non-solvable
kernel for which the existence of such a family has been established.

\end{abstract}

\section{Introduction}
\label{S.intro}

Smoluchowski's coagulation equation provides a mean-field description of
binary coalescence of clusters. 
If $\xi$ denotes the size of a cluster and $f(\xi,t)$ the corresponding number density
at time $t$  then the equation is 
\begin{equation}\label{coag1}
\frac{\partial}{\partial t} f (\xi,t) = \tfrac 1 2 \int_0^\xi \,d\eta\, K(\xi-\eta,\eta)
f(\eta,t) f(\xi{-}\eta,t)
\,-
\,  f(\xi,t) \int_0^{\infty}\,d\eta\, K(\xi,\eta) f(\eta,t)\,,
\end{equation}
where $K(\xi,\eta)$ is a kernel that describes the rate of the coalescence process.

Here we consider a specific diagonal kernel of homogeneity $\gamma < 1$, 
given by $K(\xi,\eta) = \delta_{(\xi-\eta)} \xi^{1+\gamma}$, 
that reduces (\ref{coag1}) to
\begin{equation}\label{coag2}
\frac{\partial}{\partial t} f (\xi,t) = \frac 1 4 \Big( \frac{\xi}{2},t\Big)^{1+\gamma} f^2\Big(\frac{\xi}{2},t
\Big) - \xi^{1+\gamma} f^2(\xi,t)\,.
\end{equation}

In the following we study self-similar solutions of (\ref{coag2}).
Such solutions are of the form
\begin{equation}\label{sss1}
f(\xi,t) = t^{-(1+(1+\gamma)\beta)} g\Big( \frac{\xi}{t^{\beta}}\Big)
\end{equation}
for some positive $\beta$, where
 $g$ satisfies, with $x=\xi/t^{\beta}$, that
\begin{equation}\label{sss2}
-\big( 1 + (1+\gamma)\beta \big)g - \beta x g'(x) 
= \frac 1 4 \Big( \frac{x}{2}\Big)^{1+\gamma} g^2\left( \frac{x}{2}\right)
-x^{1+\gamma} g^2(x)\,.
\end{equation}
If one looks for solutions with conserved mass, then $\beta$
is uniquely determined by $\beta=\beta_*:= 1/(1{-}\gamma)$.
For further reference we also note that we can integrate the equation in (\ref{sss2})
to obtain
\begin{equation}\label{sss2b}
\beta x^2 g(x) = \int_{x/2}^x s^{2+\gamma} g^2(s)\,ds + (1{-}\gamma)(\beta - \beta_*)
\int_0^x s g(s)\,ds\,.
\end{equation}
Here we assumed implicitly that $x g(x)$ and $x^{2+\gamma} g^2(x)$ are 
 integrable at zero and that
$\lim_{x\to 0} x^2 g(x)=0$. As we will see below (cf. (\ref{sol1})), these properties
 will be satisfied by the solutions we are going to consider.
Notice also that we have the well-known power-law solution 
\begin{equation}\label{sss5}
g = x^{-(1+\gamma)} \frac{1}{1{-}\theta} \qquad \mbox{ with } \theta:=2^{\gamma-1}<1\,.
\end{equation}

In \cite{Le06} a 
mass-conserving solution of (\ref{sss2b}), that is a solution for $\beta = \beta_*$,
is constructed that is decaying exponentially fast and satisfies
\begin{equation}\label{le1}
g(x) = x^{-(1+\gamma)} \Big( \frac{1}{1{-}\theta} - c x^{\mu/(1-\gamma)} + o(x^{\mu/(1-\gamma)})
\Big) \qquad \mbox{ as } x \to 0\,,
\end{equation}
where $\mu>0$ satisfies a certain transcendental equation. The constant $c>0$ 
is not determined due to an invariance of (\ref{sss2}) under the rescaling
$g(x) \mapsto a^{1+\gamma} g(ax)$ for any $a>0$. In the case of mass-preserving solutions
the constant can be fixed by normalizing the  mass of the solution.
As is pointed out in \cite{Le06}, the solution is unique in the class of
 functions satisfying (\ref{le1}), but uniqueness in general is not known.

In \cite{Le06} the question is raised whether solutions with algebraic decay, others from the
one in (\ref{sss5}),  exist in analogy to the ones that have been
found in \cite{MP1} for the constant and additive kernel. More precisely, for example for the constant
kernel, it is established in \cite{MP1} that there exists a family of self-similar solutions with infinite
mass and the decay behavior $x^{-(1+\rho)}$ for all $\rho \in (0,1)$. Furthermore, it is
shown that a solution
of the coagulation equation converges to the self-similar solution with decay behavior 
 $x^{-(1+\rho)}$  if and only if the
 mass-distribution of the initial data  is regularly varying with exponent $\rho$. In this note we prove
for the diagonal kernel
the existence of a corresponding family of self-similar solutions with infinite mass and asymptotic behavior  $x^{-(1+\rho)}$ as $x \to \infty$ 
with $\rho \in (\gamma,1)$. Notice, that this includes solutions that are 
increasing as $x \to \infty$ if $\gamma <-1$. 
Our proof is simple and exploits strong monotonicity properties
of a suitably rescaled version 
of the equation for the self-similar solution. We presently do not know, however, how 
to characterize the domains of attraction of these self-similar solutions. The analysis in \cite{MP1}
relies on the fact that the Laplace transform of the equation satisfies a simple ODE, a method
that is not applicable in the present situation.

Our main result is the following.
\begin{theorem}
\label{T.1}
Let $\gamma <1$ and  $\mu$ be the unique positive solution of 
\begin{equation}\label{mueq}
\frac{1+\beta \mu}{2} = \frac{1-2^{\gamma-1-\mu}}{1-2^{\gamma-1}}\,.
\end{equation}
Then there exists for any $\beta>\beta_*$ a solution $g$ of 
(\ref{sss2b}) such that
\begin{equation}\label{sol1}
g(x) = x^{-(1+\gamma)} \Big( \frac{1}{1{-}\theta} - c x^{\mu/(1{-}\gamma)} + 
o\big (x^{\mu/(1{-}\gamma)} \big)\Big) 
\end{equation}
as $x \to 0$ with a positive constant $c$ Furthermore, $x^{-(1+\gamma)} g(x)$ is
monotonically decreasing and satisfies 
\begin{equation}\label{sol2}
g(x) \sim \frac{d}{x^{1+\gamma+1/\beta}} \qquad \mbox{ as } x \to \infty
\end{equation}
for some positive constant $d$.
\end{theorem}

As explained above, the constants $c$ and $d$ in Theorem \ref{T.1} are not determined
due to the invariance of the equation under appropriate rescaling.

\section{Proof}
Our proof proceeds similarly  to the one in \cite{Le06} for the mass-conserving solutions.
First, to  scale out the singular behavior as $x \to 0$, we introduce
$h(x)=g(x) x^{1+\gamma}
$
such that
$h$ solves
\begin{equation}\label{sss4}
-\beta x h'(x) - h(x) = \theta
h^2 \left ( \frac{x}{2} \right) - h^2(x)\,
\end{equation}
or, due to (\ref{sss2b}), 
\begin{equation}\label{sss4b}
\beta x^{1-\gamma} h(x) = \int_{x/2}^x s^{-\gamma} h^2(s)\,ds + (1{-}\gamma)(\beta -\beta_*) \int_0^x s^{-\gamma}
h(s)\,ds\,.
\end{equation}
Notice, that the power-law solution (\ref{sss5}) corresponds to the constant solution
$h \equiv 1/(1{-}\theta)$. It is also clear that any solution of $(\ref{sss4})$ for which
$\lim_{x \to 0} h(x) $ exists, that this limit must equal $1/(1{-}\theta)$. We are now looking
for solutions that bifurcate from this constant at $x \to 0$. 
 
In order to identify the next order behavior, we make the ansatz $h(x) = 1/(1{-}\theta) + x^{\mu} +o(x^{\mu})$
as $x \to 0$. Plugging this into (\ref{sss4b}), 
recalling that $\beta_*=1/(1{-}\gamma)$ and rearranging we find that $\mu$ must indeed satisfy (\ref{mueq}).
If we denote by  $F(\mu)=(1-2^{\gamma-1-\mu})/(1{-}\theta)$ we see that $F(0)=1>1/2$. On the other hand,
$F$ is increasing and $\lim_{\mu \to \infty} F(\mu) = 1/(1{-}\theta)$. Hence, there must be a unique
positive solution of (\ref{mueq}).

Next, we introduce the function $j(x)$ via 
\begin{equation}\label{jdef}
h(x)=\frac{1}{1{-}\theta} + x^{\mu}\big( -c + j(x) \big)\,,
\end{equation}
where $c \in \mathbb{R}$ is a constant.
Using  equations (\ref{mueq}) and (\ref{sss4b}) we obtain that $j$ satisfies 
\begin{equation}\label{jeq}
\begin{split}
j(x)&= \frac{1}{\beta} x^{-(1-\gamma+\mu)} \Big( \int_{x/2}^x s^{-\gamma+\mu} \frac{2}{1{-}\theta}
j(s)\,ds + \int_{x/2}^x s^{-\gamma+ 2\mu} (-c+j(s))^2\,ds\\
&\quad + (1{-}\gamma)(\beta{-}\beta_*) \int_0^x s^{-\gamma+\mu} j(s)\,ds \Big) =: T[j]\,.
\end{split}
\end{equation}
In order to prove that a local solution of (\ref{jeq}) exists, we can proceed analogously to \cite{Le06}.
We only indicate the main steps here.

We define  for some $\eps \in (0,\mu)$ and $z>0$ the space
\[
C_{\eps}(z):=\{ f \in C[0,z]\,;\, f(0)=0\,;\, \|f\|:=\sup_{x \in [0,z]} x^{-\eps} |f(x)| <\infty\}\,.
\]
It is clear that the operator $T$ maps $C_{\eps}(z)$ into itself. Next, we are going to show that
$T$ maps a ball in $C_{\eps}(z)$ of a sufficiently small
radius $R$ into itself if $z$ is sufficiently small.
This follows from
\[
\begin{split}
\|T[j]\|& \leq \frac{1}{\beta}\|j\| \Big( \frac{2}{1{-}\theta} 
\frac{1}{1{-}\gamma{+}\mu{+}\eps}
\big( 1 {-} 2^{\gamma-1-\mu-\eps}\big) 
 + \|j\|  \frac{2 z^{\mu}}{1{-}\gamma{+}2\mu{+}2\eps} \big( 1 - 2^{\gamma-1-\mu-2 \eps}\big)
\\
& \qquad + c^2 \frac{ 2z^{\mu}}{1{-}\gamma{+}2\mu} + \frac{(1{-}\gamma)
(\beta{-}\beta_*) }{1{-}\gamma{+}\mu{+}\eps} \Big)
\end{split}
\]
that implies
\[
\|T[j]\| \leq \|j\| \Big(\frac{1}{\beta(1{-}\gamma{+}\mu{+}\eps)} \big( 2 F(\mu+\eps) + (1{-}\gamma) \beta -1\Big) 
+ C z^{\mu} \big( \|j\|^2 +1\big)\,.
\]
Now we know by the definition of $\mu$  that $2 F(\mu+\eps) < 1+\beta(\mu+\eps)$
and hence 
\[
\frac{1}{\beta(1{-}\gamma{+}\mu{+}\eps)} (2 F(\mu+\eps) + (1{-}\gamma) \beta -1)
< \frac{1}{1{-}\gamma{+}\mu{+}\eps }(\mu{+}\eps {+}1{-}\gamma)=1\,.
\]
Thus, there exists a constant $\kappa=\kappa(\eps)<1$ such that 
if $\|j\| \leq R$ we find
$\|T[j]\| \leq  \kappa R + Cz^{\mu}(R^2+1)$.
For sufficiently small $z$ and an appropriately small $R$ the right hand side
is bounded by $R$.  Similarly one can show that $T$ is a contraction, we omit
the details here. 
Hence, a local solution to (\ref{jeq}) exists, and thus also to (\ref{sss4}).
Next, we choose $c>0$, and claim that $h$ is decreasing in a neighbourhood of zero. 
To see this, notice that it follows from (\ref{jeq})   that $j'(x)$ exists for $x>0$ and that we have
the estimate
$|j'(x)| \leq C \frac{|j(x)|}{x} + C x^{\mu-1}$  for $x \in (0,z)$.
This in turn implies that 
\[
h'(x) = \mu x^{\mu-1} \big( -c+j(x)) + x^{\mu} j'(x) \leq x^{\mu-1} \big( - c \mu + C \mu|j(x)| + C x^{\mu}
\big)\,.
\]
If $z$ is sufficiently small, we find that $h'(x)<0$ for $x \in (0,z)$.
We are going to show that as long as $h$ exists and is positive this property is
conserved. Indeed, assume that there exists
 $x_0>0$  such that $h'(x_0)=0$. Then (\ref{sss4}) and the fact
that $h$ is decreasing for $x<x_0$ imply
that 
\[
0=  h(x_0)^2 - h(x_0)- \theta h^2\Big( \frac{x_0}{2}\Big)^2
< \big(1{-}\theta\big) h^2(x_0) - h(x_0)
=h(x_0) \big( (1{-}\theta) h(x_0) -1\big) \,.
\]
As long as $h$ is positive, the right hand side is strictly negative, since 
$h(x_0)<1/(1{-}\theta)$ and we obtain the desired contradiction.
Moreover,  equation (\ref{sss4b}) implies for $\beta \geq \beta_*$
that $h$ is positive whenever it exists. Hence, using standard results on ordinary differential equations,
we obtain global existence  of a solution $h$
to (\ref{sss4}) which is strictly decreasing.
Since equation (\ref{sss4}) has the only stationary points $1/(1{-}\theta)$ and $0$, it also follows
that $h(x) \to 0$ as $x \to \infty$. 

It remains to show that $h(x) \sim  d x^{-1/\beta}$ as $x \to \infty$ from which (\ref{sol2}) follows. First, due to the invariance of equation (\ref{sss4})
under the transformation $x \to ax$ for $a>0$, we can assume without
loss of generality that $h(1)=1/2$. Since
$h$ satisfies
$\beta x h'(x) + h(x) \leq h^2(x)$
we have by simple comparison that
\begin{equation}
\label{supersol}
h(x) \leq 
\frac{1}{1+  x^{1/\beta}} \qquad \mbox{ for } x\geq 1\,.
\end{equation}
We now  introduce $p(x)=x^{1/\beta}h(x)$ that solves
\begin{equation}\label{peq}
\beta p'(x) = x^{-(1+1/\beta)} \Big( p^2(x) - \theta 2^{2/\beta} p^2
\Big (\frac{x}{2}\Big) \Big )\,.
\end{equation}
The estimate (\ref{supersol}) in particular implies that $p(x) \leq 1$
for all $x \geq 1$ and thus (\ref{peq}) implies that $\beta |p'(x)|
\leq 2 x^{-(1+1/\beta)} $ for all $x \geq 2$. Hence
$|p(x)-p(x_0)| \leq 2 x_0^{-1/\beta}
$
for any $x_0\geq 2$ which implies that $\lim_{x\to \infty} p(x)$ exists. In order
to complete the proof of Theorem \ref{T.1} it remains
to establish that this limit is strictly positive. To this end we note that
(\ref{sss4b}) implies 
\begin{equation}\label{hineq}
\beta x^{1-\gamma }h\left( x\right) >\left( 1{-}\gamma \right) \left( \beta
-\beta _{\ast }\right) \int_{0}^{x}s^{-\gamma }h\left( s\right) ds\,.
\end{equation}
If we define
$\Phi \left( x\right):  =\int_{0}^{x}s^{-\gamma }h\left( s\right) ds$ then (\ref{hineq})
implies that $
\beta x\Phi ^{\prime }\left( x\right)$\linebreak $ -\left( 1-\gamma \right) \left( \beta
-\beta _{\ast }\right) \Phi \left( x\right) >0$. Integrating this last inequality we obtain \linebreak
$\left( x^{-\frac{\left( 1-\gamma \right) \left( \beta -\beta _{\ast }\right) 
}{\beta }}\Phi \left( x\right) \right) ^{\prime }>0
$ and  thus
\[
x^{-\frac{\left( 1-\gamma \right) \left( \beta -\beta _{\ast }\right) }{%
\beta }}\Phi \left( x\right) \geq \Phi(1) = \int_0^1 s^{-\gamma} h(s)\,dx=:c_0 >0\,
\]%
for all $x\geq 1$.
Thus
\[
\Phi \left( x\right) \geq c_0 x^{\frac{\left( 1-\gamma \right) \left( \beta
-\beta _{\ast }\right) }{\beta }}=c_0 x^{1-\gamma }x^{-\frac{1}{\beta }}\ \ \
\]%
for $x \geq 1$ and plugging this into (\ref{hineq}) we find
$h\left( x\right) \geq \frac{c_0}{\beta }x^{-\frac{1}{\beta }}$
for all $x \geq 1$, that finishes the proof.

\bigskip
{\bf Acknowledgment.} This work was supported by the EPSRC Science and Innovation award to the
Oxford Centre for Nonlinear PDE (EP/E035027/1),
the DGES Grant MTM2007-61755, the Proyecto Intramural 200850I248
 and the Isaac Newton Institute.

\bibliographystyle{amsplain}

\end{document}